\documentclass[11pt]{amsart}
\usepackage{amssymb}

\usepackage{palatino}

\newtheorem{thm}{Theorem}
\newtheorem{dfn}{Definition}
\newtheorem{rmk}{Remark}

\newtheorem{prop}{Proposition}
\newtheorem{lem}{Lemma}
\newtheorem{eg}{Example}

\newcommand{\Lg}{\mbox{$\mathfrak{g}$}}
\newcommand{\Lk}{\mbox{$\mathfrak{k}$}}
\newcommand{\Lp}{\mbox{$\mathfrak{p}$}}

\newcommand{\La}{\mbox{$\mathfrak{a}$}}

\newcommand{\Lb}{\mbox{$\mathfrak{b}$}}
\newcommand{\Lh}{\mbox{$\mathfrak{h}$}}

\newcommand{\Lc}{\mbox{$\mathfrak{c}$}}
\newcommand{\Ll}{\mbox{$\mathfrak{l}$}}

\newcommand{\Lv}{\mbox{$\mathfrak{v}$}}
\newcommand{\Lt}{\mbox{$\mathfrak{t}$}}
\newcommand{\Ad}{\mbox{Ad}}

\newcommand{\Pf}{{\em Proof}. }
\newcommand{\EPf}{\hfill$\square$}

\newcommand{\R}{\mbox{$\mathbf R$}}
\newcommand{\C}{\mbox{$\mathbf C$}}

\newcommand{\ontop}[2]{\genfrac{}{}{0pt}{}{#1}{#2}}

\begin{document}


\title{The discriminants associated to\\isotropy representations
of symmetric spaces}


\author{Claudio Gorodski}\thanks{This work has been inspired by a talk given by 
Peter Lax at the Institute of Mathematics and Statistics of the 
University of S\~ao Paulo in April 2010.}

\date{\today} 

\address{Instituto de Matem\' atica e Estat\'\i stica\\
         Universidade de S\~ ao Paulo\\
         Rua do Mat\~ao, 1010\\
         S\~ ao Paulo, SP 05508-090\\
         Brazil}
\email{\texttt{gorodski@ime.usp.br}}

\begin{abstract}
We consider a generalized discriminant associated to a symmetric 
space which generalizes the discriminant of real symmetric matrices,
and note that it can be written as a sum of squares of real polynomials. 
A method to estimate the minimum number of squares required to 
represent the discrimininant is developed and applied in examples.    
\end{abstract}

\maketitle 


\section{Introduction}

In his doctoral dissertation held in K\"onigsberg 
in 1885, Minkowski proposed the conjecture that, unlike the quadratic case, 
nonnegative homogeneous polynomials 
of higher degree and more than two variables in general cannot 
be written as a sum of squares of real polynomials. The problem attracted
the attention of Hilbert who in 1888 proved nonconstructively 
the existence of such 
polynomials. However, the first concrete example of a nonnegative polynomial
which is not a sum of squares  
seems to have been given only in 1967 by Motzkin~\cite{mot}. 
The question of which nonnegative polynomials admit such representations
is of interest in real algebraic geometry 
and practical importance in applied mathematics,
but it is left open. On the other hand, 
Blekherman~\cite{ble} has shown that  
there are significantly more nonnegative polynomials than sums of squares
of polynomials by computing asymptotic bounds on the sizes 
of these sets. 

The discriminant of a real symmetric $n\times n$ matrix $Y$ is 
\[ \delta(Y) = \Pi_{i<j}(\theta_i(Y)-\theta_j(Y))^2, \]
where the $\theta_i(Y)$ are the eigenvalues of $Y$. It is well 
known that the discriminant
$\delta(Y)$ is a nonnegative homogeneous polynomial of degree $n(n-1)$
in the entries of $Y$, and it vanishes if and only if $Y$ has an
eigenvalue of multiplicity bigger than one. 
In the nineteenth century, Kummer exhibited 
an explicit representation of $\delta$ as a sum of squares for 
$n=3$ (the case $n=2$ is immediate), 
and Borchardt generalized it for arbitrary $n$. 
More recently, several authors have rediscovered and refined
these results in one or another form 
(see~\cite{wat,new,ily,lax,par,dom} and the references 
therein). In particular,
the approaches of Lax~\cite{lax} and Domokos~\cite{dom}
(see also~\cite{wat})
make important use of the action of the orthogonal group
on the space of symmetric matrices by conjugation:
since conjugate matrices have the same set of eigenvalues, 
this action leaves the discriminant unchanged.

In this note, we remark that some of these results can be viewed
in the realm of symmetric spaces (or, slightly more generally,
polar representations). Indeed the isotropy representations
of Riemannian symmetric spaces constitute a remarkable class of 
representations of compact Lie groups (see e.g.~\cite[ch.~3]{bco} for
a discussion). Herein we are interested in the functional which computes
the volume of the orbits of those representations; it turns out
that its square is an invariant polynomial that can be considered 
as a generalized discriminant associated to the symmetric space 
(the case of symmetric matrices alluded to above corresponds to 
the symmetric space $GL(n,\R)/SO(n)$, see example below).

More generally, the following result is probably well known, 
but perhaps has not been related to the present context.

\begin{thm}\label{disc}
The functional computing the squared volume
of the principal orbits of an orthogonal 
representation of a compact connected Lie group uniquely extends 
to an invariant homogeneous polynomial function on the representation 
space. Moreover, this polynomial can be expressed as a sum of squares
of polynomials.
\end{thm}

The number of squares involved in this representation 
is big (for instance, for isotropy representations
of symmetric spaces of maximal rank this number
is the binomial coefficient $\binom dr$, where $d$ and $r$ 
are respectively the dimension and rank of the symmetric space). 
In the case of polar representations, it is possible
to show the existence of expressions with much smaller
number of squares, based on Theorem~\ref{decomp} below.
In particular, for isotropy representations
of symmetric spaces of maximal rank, we prove Theorem~\ref{alg} below, 
which supplies an effective method to construct those
expressions. Theorems~\ref{decomp} and~\ref{alg} 
generalize some ideas from~\cite{lax,dom}. By applying this 
method, we recover

\begin{thm}[\cite{dom}]\label{dom} 
The discriminant of $n\times n$ real symmetric matrices (with 
respect to the orthogonal group) can be written as the sum of 
$\binom{2n-1}{n-1}-\binom{2n-3}{n-1}$ squares.
\end{thm}

We also apply our method to complex symmetric matrices, 
or the symmetric space $Sp(n,\R)/U(n)$. We prove

\begin{thm}\label{spn-un}
The discriminant of $n\times n$ complex symmetric matrices
(with respect to the unitary group) can be written as the sum
of $2\,\binom{2n-1}{n}$ squares.
\end{thm}

\emph{Conventions.} Let a compact connected 
Lie group $G$ act linearly on a real vector space $V$.
We consider the space $\R[V]$ of real polynomials
on $V$, and note that $G$ acts on a polynomial $f$ by the rule
$(gf)(v)=f(g^{-1}v)$, where $v\in V$. The space
of invariants is denoted $\R[V]^G$.
The homogeneous component of degree $m$ of $\R[V]$ is denoted $\R[V]_m$.
An element $v\in V$ 
and its orbit $Gv$ are called
\emph{regular} if $\dim Gv$ is maximal amongst $G$-orbits, and 
\emph{singular} otherwise. Suppose $v$ is a regular element; then 
the dimension of its isotropy group $G_v$ is as small as possible.
In this case, $Gv$ is called a \emph{principal} orbit if in addition the 
number of connected components of $G_v$ is as small as possible;
otherwise, $Gv$ is called an \emph{exceptional} orbit. 
We also consider the complexified representation 
on $V^c$ (of $G$ or $G^c$) and the space of complex polynomials 
$\C[V^c]$. On the other hand, if $W$ is a already a 
complex representation, $W^r$ denotes its realification. 

\section{The volume functional as a sum of squares} 

Let $V$ be a real orthogonal finite-dimensional
representation space of a compact connected
Lie group $G$. The principal orbit type theorem asserts that the 
union of all principal orbits is an open, dense, 
invariant set with connected
orbit space. Indeed every principal orbit is a homogeneous Riemannian
manifold of $G/H$-type with $G$-invariant metric induced from $V$,
where $H$ is a fixed closed subgroup of $G$
(more generally,
if $v$ is a regular point, the orbit $Gv$ is finitely covered
by $G/H$), 
and every orbit meets the fixed point subspace $V^H$.  
Let $\Lg$ be the Lie algebra of $G$ equipped with
an $\mathrm{Ad}_H$-invariant inner product, and let $\Lh$ be the 
Lie algebra of $H$. Fix orthonormal bases $x_1,\ldots,x_m$ of $\Lh^\perp$
and $v_1,\ldots,v_d$ of $V$. Denote by $\tilde x_1,\ldots,\tilde x_m\in
\mathfrak{so}(V)$ the induced Killing fields on $V$, namely 
$\tilde x(v)=x_i\cdot v$ for $v\in V$. For a point $v\in V^H$,
it is easy to compute the Jacobian determinant of the $G$-equivariant
diffeomorhism $G/G_v\to Gv$ to deduce that the
$m$-dimensional volume of the orbit $Gv$ is
\[ \mathrm{vol}(Gv)=\left\{\begin{array}{cl}
||\tilde x_1(v)\wedge\cdots\wedge\tilde x_m(v)||\,\mathrm{vol}(G/H)
&\quad\mbox{if $Gv$ is a principal orbit,}\\
\frac 1k||\tilde x_1(v)\wedge\cdots\wedge\tilde x_m(v)||\,\mathrm{vol}(G/H)&\quad\mbox{if $Gv$ is an exceptional orbit,}\\
0&\quad\mbox{if $Gv$ is a singular orbit,}
\end{array}\right. \]
where $k$ is the index of the covering $G/H\to Gv$.  
The value of the constant $\mathrm{vol}(G/H)$ is unimportant, and in  
the sequel we shall normalize the metric in $\Lg$ so that this constant
becomes $1$. 

\begin{dfn}
\em
The \emph{discriminant} of the representation
$(G,V)$ is the $G$-invariant function on $V$ defined by
\[ \delta(v) = \left\{\begin{array}{cl}
k^2\,\mathrm{vol}^2(Gv)&\quad\mbox{if $v$ is a regular point,}\\
0&\quad\mbox{if $v$ is a singular point,}
\end{array}\right. \]
where $v\in V$ and $k$ is the index of the covering $G/H\to Gv$
($k=1$ if $Gv$ is a principal orbit;
in particular, if there are no exceptional orbits, then 
$\delta(v)=\mathrm{vol}^2(Gv)$ for all $v\in V$). 
\end{dfn}

If $v\in V^H$, the formula
\[ \delta(v) = ||\tilde x_1(v)\wedge\cdots\wedge\tilde x_m(v)||^2
=\det\left(\langle\tilde x_i(v),\tilde x_j(v)\rangle\right) \]
shows that $\delta|_{V^H}$ is a
homogeneous polynomial of degree $2m$ on $V^H$. The restriction
$\delta|_{V^H}$ must also be invariant under the subgroup of $G$
stabilizing $V^H$, which is the normalizer $N(H)$ of $H$ in $G$. 
The Luna-Richardson theorem~\cite{lr}, 
asserting that the restriction map from
$V$ to $V^H$ induces a graded algebra isomorphism $\R[V]^G\to\R[V^H]^{N(H)}$,
can now be used to conclude that $\delta$ is a 
homogeneous polynomial of degree $2m$ on $V$.
Recall that by a result of Hilbert and Hurwitz, the  
algebra $\R[V]^G$ is finitely generated.  

On the other hand, it is very easy to make
explicit the polynomial structure of $\delta$. 
We shall recall that while introducing some new objects. 
Take an orthonormal basis $\{v_{i_1}\wedge\cdots\wedge v_{i_m}\}$ 
of $\Lambda^mV$,
where $\{i_1<\cdots<i_m\}$ is an increasing multi-index, and write,
for $v\in V^H$,
\begin{eqnarray*}
\delta(v) & = & ||\tilde x_1(v)\cdots\wedge\tilde x_m(v)||^2 \\
          & = & \sum_{1\leq i_1<\cdots<i_m\leq d}
\langle \tilde x_1(v)\wedge\cdots\wedge\tilde x_m(v),v_{i_1}\wedge\cdots\wedge v_{i_m}
\rangle^2 \\
&=& \sum_{1\leq i_1<\cdots<i_m\leq d}\left[\det||\,\langle x_b\cdot v,v_{i_a}\rangle\,||_{a,b=1}^m\right]^2.
\end{eqnarray*}

Next, we consider the morphism $\rho:V\to\mathrm{Hom}(\Lg,V)$ encoding
the representation, namely $\rho(v)(x)=x\cdot v$ for $v\in V$, $x\in\Lg$. 
Plainly, $\rho$ is $G$-equivariant: 
$\rho(gv)=g\rho(v)\mathrm{Ad}_{g^{-1}}$ for $g\in G$. 
Extend $\{x_j\}_{j=1}^m$ to a basis $\{x_j\}_{j=1}^p$ of $\Lg$. 
Then $\rho(v)$ is represented by the matrix $(\rho(v)_{ij})$, where
$\rho(v)(x_j)=\sum_{i=1}^d\rho(v)_{ij}v_i$ for $j=1,\ldots,p$. 
Since $\langle x_b\cdot v,v_{i_a}\rangle=\rho(v)_{i_a,b}$, we have
\[ \delta(v)= \sum_{1\leq i_1<\cdots<i_m\leq d}\left[\det||\,\rho(v)_{i_a,b}\,||_{a,b=1}^m\right]^2 \] 
for $v\in V^H$ and, in fact, we can write 
\begin{equation}\label{delta}
\delta(v)= \sum_{\ontop{1\leq i_1<\cdots<i_m\leq d}{1\leq j_1<\cdots<j_m\leq d}}\left[\det||\,\rho(v)_{i_a,j_b}\,||_{a,b=1}^m\right]^2
\end{equation}
as the last $p-m$ columns of $(\rho(v)_{ij})$ contain only
zeros. 

Define a morphism $\Phi:\Lambda^m\Lg\otimes\Lambda^mV\to\R[V]_m$
by setting 
\[ \Phi(x_{j_1}\wedge\cdots\wedge x_{j_m}\otimes v_{i_1}\wedge\cdots\wedge v_{i_m})(v)=\langle x_{j_1}(v)\wedge\cdots\wedge x_{j_m}(v),v_{i_1}\wedge\cdots\wedge
 v_{i_m}\rangle. \]
Then $\Phi$ is equivariant, so by Schur's lemma, it maps a $G$-irreducible
component of $\Lambda^m\Lg\otimes\Lambda^mV$ either to zero or onto an
isomorphic $G$-irreducible component of $\R[V]_m$. 
Note that
$\Phi(x_{j_1}\wedge\cdots\wedge x_{j_m}\otimes v_{i_1}\wedge\cdots\wedge v_{i_m})(v)=\det||\rho(v)_{i_a,j_b}||_{a,b}$ for all $v\in V$, so
the image of $\Phi$ lies in the ideal of polynomials vanishing
along the variety of singular points. 

The following lemma has been used in different contexts 
(e.g.~\cite[Lemma~2.1]{dom} or~\cite[p.~53]{bco}).

\begin{lem}\label{inv}
Let $V$, $W$ be real orthogonal representations of a compact Lie
group $G$, and let $\Psi:W\to\R[V]$ be a $G$-equivariant 
map. Then, for an orthonormal basis $\{w_i\}$ of $W$, the 
polynomial $f=\sum_i\Psi(w_i)^2$ is $G$-invariant and independent
of the basis used to construct it.
\end{lem}

\Pf Let $v\in V$, $g\in G$. Since the $G$-action on $W$ is orthogonal,
there is an orthogonal matrix $(a_{ij})$ such that $gw_j=\sum_i a_{ij}w_i$.
Now
\begin{eqnarray*}
 f(g^{-1}v) & = & \sum_j\Psi(w_j)(g^{-1}v)^2 \\
           & = & \sum_j\Psi(gw_j)(v)^2 \\
           & = & \sum_j\left[\sum_ia_{ij}\Psi(w_i)(v)\right]^2 \\
           & = & \sum_{i,k}\left(\sum_ja_{ij}a_{kj}\right)\Psi(w_i)(v)\Psi(w_k)(v) \\
& = & \sum_i \Psi(w_i)(v)^2 \\
&= & f(v),
\end{eqnarray*}
as desired. The last assertion in the statement is proven similarly. \EPf

\medskip

It follows from the lemma that the right hand-side
of equation~(\ref{delta}) represents a $G$-invariant polynomial
on $V$. Hence we have found a polynomial expression for $\delta$, 
namely
\[  \delta = \sum_{I,J}\left[\Phi(x_J\otimes v_I)\right]^2, \] 
where $I$, $J$ are increasing multi-indices; indeed, this is a sum of squares
whose number is the dimension of the domain $\Lambda^m\Lg\otimes\Lambda^mV$ 
of $\Phi$. This proves Theorem~\ref{disc} stated in the introduction. 
The lemma also says that,
for every $G$-irreducible subspace $W$ of $\Lambda^m\Lg\otimes\Lambda^mV$ 
with $\Phi(W)\neq 0$, and $\{w_i\}$ an orthonormal basis of $W$, 
the polynomial $f_W:=\sum_i\Phi(w_i)^2$ is $G$-invariant. 

\subsection{The case of polar representations}

We can get better results if we assume that $(G,V)$ is polar,
as we henceforth do. This means there exists a subspace $\Lc\subset V$,
called a Cartan subspace, that meets every $G$-orbit, and meets 
always orthogonally. We can choose the basis of $V$ so that 
$v_1,\ldots,v_r\in\mathfrak c$, where $r=\dim\mathfrak c$. Then the tangent
spaces $T_v(Gv)$ for regular
$v\in\mathfrak c$ are all parallel and indeed
spanned by $v_{r+1},\ldots, v_d$. 
Now the matrix $\rho(v)$ for regular
$v\in\Lc$ has the block form 
$\left(\begin{array}{cc}0&0\\ \star&0\end{array}\right)$.
We recall that polar representations
do not admit exceptional orbits~\cite[Cor.~5.4.3]{bco}. 

Consider the special element 
\[ \vartheta=x_1\wedge\cdots\wedge x_m\otimes
v_{r+1}\wedge\cdots\wedge v_d\in
\Lambda^m\Lg\otimes\Lambda^mV. \]
It follows from the block form of the above matrix
that the restriction map $V\to\Lc$ takes
all determinants of minors $\Phi(x_J\otimes v_I)$ to zero but
$\Phi(\vartheta)$. Therefore
\[ \delta|_{\mathfrak c}=\Phi(\vartheta)^2|_{\mathfrak c} \]
and 
\begin{eqnarray*}
f_W|_{\mathfrak c} & = & \sum_i[\Phi(w_i)|_{\mathfrak c}]^2 \\
                & = & \sum_i[\langle w_i,\vartheta\rangle
\Phi(\vartheta)|_{\mathfrak c}]^2 \\
                & = & \left(\sum_i\langle w_i,\vartheta\rangle^2
\right)[\Phi(\vartheta)|_{\mathfrak c}]^2 \\
                & = & c\cdot \delta|_{\mathfrak c}, 
\end{eqnarray*}
where $c$ denotes a positive constant. 
Since $f_W$ and $\delta$ are both $G$-invariant polynomials, this 
implies that $f_W=c\cdot\delta$ on $V$. 
Hence:

\begin{thm}\label{decomp}
If $(G,V)$ is polar, then 
each $G$-irreducible component of the image of $\Phi$ 
gives rise to a decompositions of $\delta$ as a sum of squares 
of polynomials.
\end{thm}

\begin{rmk}
\em
Polarity is necessary in the statement of Theorem~\ref{decomp}.
In fact, the diagonal action of $SO(2)$ on $\R^2\oplus\R^2$ is not polar.
Let $\mathfrak{so}(2)=\langle x\rangle$. Then
\[ x\cdot\left(\begin{array}{cc}a_1&b_1\\a_2&b_2\end{array}\right)
= \left(\begin{array}{cc}-a_2&-b_2\\a_1&b_1\end{array}\right) \]
and
\[ \rho\left(\begin{array}{cc}a_1&b_1\\a_2&b_2\end{array}\right)
=\left(\begin{array}{c}-a_2\\a_1\\-b_2\\b_1\end{array}\right). \]
Of course, $\delta=a_1^2+a_2^2+b_1^2+b_2^2$ is $SO(2)$-invariant. 
On the other hand, $m=1$ and $\Lambda^m\Lg\otimes\Lambda^mV=
\mathfrak{so}(2)\otimes(\R^2\oplus\R^2)\cong\R^2\oplus\R^2$. 
By taking $W=\R^2\oplus0$ or $0\oplus\R^2$, we get the 
$SO(2)$-invariant $f_W=a_1^2+a_2^2$ or $b_1^2+b_2^2$, each of
which is different from $\delta$.
\end{rmk}

\medskip\noindent\textbf{Baby example.}
The direct product action $(SO(2)\times SO(2),\R^2\oplus\R^2)$ is 
clearly polar with rank $r=2$.
Let $x_1$, $x_2$ be generators of the summands
in $\mathfrak{so}(2)+\mathfrak{so}(2)$, let $\{v_1,v_2;v_3,v_4\}$
be an orthonormal basis of $\R^2\oplus\R^2$, and let 
$\{a_1,a_2;b_1,b_2\}$ be the dual basis. Then
\[ x_1\cdot\left(\begin{array}{cc}a_1&b_1\\a_2&b_2\end{array}\right)
= \left(\begin{array}{cc}-a_2&0\\a_1&0\end{array}\right)\quad\mbox{and}\quad
x_2\cdot\left(\begin{array}{cc}a_1&b_1\\a_2&b_2\end{array}\right)
= \left(\begin{array}{cc}0&-b_2\\0&b_1\end{array}\right), \]
and
\[ \rho\left(\begin{array}{cc}a_1&b_1\\a_2&b_2\end{array}\right)
=\left(\begin{array}{cc}-a_2&0\\a_1&0\\0&-b_2\\0&b_1\end{array}\right). \]
Here $\delta=a_1^2b_1^2+a_1^2b_2^2+a_2^2b_1^2+a_2^2b_2^2=(a_1^2+a_2^2)(b_1^2+b_2^2)
=f_1f_2$, where $f_1=a_1^2+a_2^2$, 
$f_2=b_1^2+b_2^2$ is a complete system of invariants. Moreover,
$m=2$ and
\[ \Lambda^2\Lg\otimes\Lambda^2V\cong\Lambda^2(\R^2\oplus\R^2)\cong
\Lambda^2\R^2\oplus\Lambda^2\R^2\oplus\R^2\otimes\R^2. \]
The first two summands on the left hand-side are spanned
by $\{v_1\wedge v_2,v_3\wedge v_4\}$ and mapped to zero under $\Phi$.
On the other hand, $\R^2\otimes\R^2$ decomposes into $\R^2\oplus\R^2$ yielding
\[ \delta = (a_1b_1+a_2b_2)^2+(a_1b_2-a_2b_1)^2=
(a_1b_1-a_2b_2)^2+(a_1b_2+a_2b_1)^2. \]

\section{Systems of restricted roots for symmetric spaces}

Dadok~\cite{dadok} has shown that a polar
representation of a compact Lie group has the same orbits as the isotropy
representation of a symmetric space, so we shall 
assume, with no loss of generality, 
that $(G,V)$ is already the isotropy representation of a 
symmetric space. This means that $\Ll:=\Lg+V$ has the structure
of a Lie algebra of which $\Lg+V$ is the decomposition into $\pm1$-eigenspaces
of an involutive automorphism of $\Ll$, and then the associated 
symmetric space is $L/G$ where $L$ is a Lie group
with Lie algebra $\Ll$. 
However, for the sake of tradition, henceforth we change the notation and write 
$(K,\Lp)$ ($K$ connected) for the basic representation, 
$\Lg=\Lk+\Lp$ for the corresponding Lie algebra with involution,
and $G/K$ for the associated symmetric space.
Now $(K,\Lp)$ is equivalent to the isotropy representation of $K$ on 
$T_{1K}(G/K)$. Since a symmetric space and its dual have 
equivalent isotropy representations, we may restrict our attention to
symmetric spaces of noncompact type. 
In this case the inner product on $\Lp$ is the restriction of the 
Cartan-Killing form of $\Lg$ to $\Lp$, and the 
inner product on $\Lk$ is the negative of the restriction of the 
Cartan-Killing form of $\Lg$ to $\Lk$. 

It is instructive to express the discriminant associated to a 
symmetric space in terms of its system of restricted roots. 
A Cartan subspace of $\Lp$ is the same as a maximal Abelian 
subspace $\La$ of $\Lp$. The system of restricted
roots of $\Lg$ with respect to $\La$ is the set of linear functionals
$\Lambda\subset\La^*\setminus\{0\}$ such that 
$\lambda\in\Lambda$ if and only if 
$\Lg_\lambda\neq0$, where 
$\Lg_\lambda=\{ x \in \Lg\;|\;\mathrm{ad}_ax=\lambda(a)x\ \mbox{for 
all $a\in\La$}\}$~\cite[ch.6, \S1]{loo}.   
We have the real orthogonal restricted root decomposition
$\Lg = \Lk_0 + \La + \sum_{\lambda\in\Lambda} \Lg_\lambda$
where $\Lk_0$ denotes the centralizer of $\La$ in~$\Lg$. Introduce a
lexicographic order in the dual $\La^*$ (with respect to some 
basis) and let $\Lambda^+$ denote the set of positive 
restricted roots. Since the involution of $\Lg$ interchanges $\Lg_\lambda$
and $\Lg_{-\lambda}$, we have 
decompositions
\[ \Lk = \Lk_0 + \sum_{\lambda\in\Lambda^+} \Lg_\lambda \quad\mbox{and}
\quad \Lp = \La + \sum_{\lambda\in\Lambda^+} \Lp_\lambda, \]
where $\Lk_\lambda=\Lk\cap(\Lg_\lambda+\Lg_{-\lambda})$,
$\Lp_\lambda=\Lp\cap(\Lg_\lambda+\Lg_{-\lambda})$, and 
$\dim\Lg_\lambda=\dim V_\lambda=m_\lambda$ is defined to be
the multiplicity of $\lambda$. 
For each $\lambda\in\Lambda^+$, 
there exist bases $\{x_{\lambda,j}\}_{j=1}^{m_\lambda}$
and $\{y_{\lambda,j}\}_{j=1}^{m_\lambda}$ of 
$\Lk_\lambda$ and $\Lp_\lambda$, resp., 
such that $\mathrm{ad}_a(x_{\lambda,j})=\lambda(a)y_{\lambda,j}$ 
and $\mathrm{ad}_a(y_{\lambda,j})=\lambda(a)x_{\lambda,j}$
for $a\in\La$~\cite[p.61]{loo}. 
It easily follows that 
\begin{equation}\label{delta-roots}
 \delta(a) = \Pi_{\lambda\in\Lambda^+}\lambda(a)^{2m_\lambda} 
\end{equation}
for $a\in\La$. 

\begin{eg}
\em
The symmetric space $GL(n,\R)/SO(n)$ has associated 
Cartan decomposition $\mathfrak{gl}(n,\R)=\mathfrak{so}(n)+\mathrm{Sym}(n,\R)$
where $\mathrm{Sym}(n,\R)$ is the space of real symmetric 
$n\times n$ matrices. A maximal Abelian subspace $\La$ is given by the 
subspace of diagonal matrices and then the restricted roots are 
$\theta_i-\theta_j$, $i\neq j$,
where $\theta_i\in\La^*$ is the $i$th diagonal coordinate,
and one can take the positive ones to correspond to $i<j$. 
All the multiplicities are one. The diagonal elements of a symmetric
matrix are its eigenvalues.  
Hence in this case the discriminant of 
$(SO(n),\mathrm{Sym}(n,\R))$ 
\[ \delta=\Pi_{i<j}(\theta_i-\theta_j)^2 \]
coincides with the usual discriminant of symmetric matrices,
where $\theta_1(Y),\ldots,\theta_n(Y)$ denote the eigenvalues of 
$Y\in\mathrm{Sym}(n,\R)$. 
\end{eg}

Some other examples of volume functionals of isotropy
representations of symmetric spaces given
in terms of formula~(\ref{delta-roots})
are listed in the final tables
of~\cite{hl} (the spaces listed there are in fact the compact dual symmetric spaces, but the corresponding isotropy representations are equivalent). 

The Weyl group of $\Lg$ with respect to $\La$, denoted $W(\La)$, 
is a finite group
generated by reflections on the singular hyperplanes in~$\La$ (which are the
kernels of the restricted roots) and can be realized
as the quotient of the normalizer of $\La$ in $K$ by the 
centralizer of $\La$ in $K$; it
acts on $\La^*$ by permuting the 
restricted roots. The Chevalley restriction theorem says that the 
restriction map from $\Lp$ to $\La$ induces an algebra isomorphism 
$\R[\Lp]^K\to\R[\La]^{W(\mathfrak a)}$. We use this result in the 
proof of the following proposition.

\begin{prop}
The discriminant of $(K,\Lp)$
is an irreducible polynomial  
if and only if the diagram of restricted roots 
is connected and does not contain a multiple link, namely it is of type
$A_n$ ($n\geq1$), $D_n$ ($n\geq4$), $E_6$, $E_7$ or $E_8$, and all the 
multiplicities are $1$.  
\end{prop}

\Pf It is clear that $\delta$ can be irreducible only if the 
representation $(K,\Lp)$ is irreducible, which is to say that
the diagram is connected; so we may as well make this assumption
throughout.  

We first claim that $\delta$ is reducible if and only if 
there exists a nontrivial factorization of $\delta|\La$ into 
$W(\La)$-invariant polynomials. In fact, 
suppose $\delta = f\cdot g$ is a nontrivial
factorization. Then $\delta=\Ad_kf\cdot\Ad_kg$ for $k\in K$ is a 
continuous family of nontrivial factorizations. Since a real polynomial
has a unique decomposition into irreducible factors, up to permutation
of the factors and multiplication by units, it is clear that 
$f$ and $g$ must be $K$-invariant polinomials (recall that $K$ is 
assumed connected); thus, $\delta|\La=(f|\La)(g|\La)$ where 
$f|\La$ and $g|\La$ are $W(\La)$-invariant. The converse 
follows from Chevalley's theorem.

On the other hand, it is known that the action of $W(\La)$ 
preserves the length and the multiplicities of the roots, and is transitive 
on the sets of roots of the same length. Moreover,
for each $\lambda\in\Lambda^+$, $W(\La)$ contains an element 
that maps $\lambda$ to $-\lambda$ and induces a permutation
on $\Lambda^+\setminus\{\lambda\}$. Now it is clear that 
a nontrivial $W(\La)$-invariant factor of 
$\delta|\La=\Pi_{\lambda\in\Lambda^+}\lambda^{2m_\lambda}$
exists if and only if $W(\La)$ is not transitive on $\Lambda$ 
or some multiplicity is bigger than $1$. We finish the proof
by noting that $W(\La)$ is transitive on $\Lambda$ precisely
if the diagram is one of those listed in the statement. \EPf

\medskip

The zero set $Z(\delta)$ of the discriminant is the singular set
$\mathcal S$, which is the union of the singular orbits,
and we now have a good description of its structure. 
The intersection $\mathcal S\cap\La$, being the the union of 
singular hyperplanes, is stratified by the intersections of the various 
subfamilies of singular hyperplanes, and 
$\mathcal S= \cup_{\lambda\in\Lambda^+}
K(\ker\lambda)$ has a natural, induced invariant 
stratification, which coincides with the stratification by orbit types.  
Denote by $(\ker\lambda)^0$ 
the open dense subset of the singular hyperplane $\ker\lambda\subset\La$
consisting of points not lying in any other singular hyperplane.
Then $K((\ker\lambda)^0)$ is a stratum of $\mathcal S$
of dimension $\dim\Lp-1-m_\lambda$, as is easy to see.
Therefore $\mathcal S$ has codimension
$1+\min_\lambda m_\lambda\geq2$. It is interesting to note that 
despite being defined as the zero set of a single polynomial, 
$\mathcal S$ has codimension bigger than one. Of course, 
this is related to the fact that 
$\delta$ is a sum of squares of real polynomials.

\section{Symmetric spaces of maximal rank} 

Although being a sum of squares is ingrained in the definition of 
$\delta$, in this section we propose to 
use Theorem~\ref{decomp} to find representations of it as  
sums of squares with as few terms as possible. 
In the case of the discriminant of real symmetric matrices,
this problem has its roots in classical papers
of Kummer and Borchardt, and more recently has been 
studied by Lax~\cite{lax} and 
Domokos~\cite{dom}. In practice, we need to decompose
$\Phi:\Lambda^m\Lk\otimes\Lambda^m\Lp\to\R[\Lp]_m$ into 
$K$-irreducible components and identify components which  
are not mapped to zero under $\Phi$. In the sequel,
we shall discuss this problem in the special case of locally free actions,
namely, when the principal isotropy is discrete. 
For the symmetric space $G/K$, this means
that $\mathrm{rank}(G/K)=\mathrm{rank}(G)$, i.e.~it is of maximal rank;
equivalently, $G/K$ has uniform multiplicities $1$. 
Then $\Lk_0=0$, $\dim\Lk=m$ and $\Phi:\Lambda^m\Lp\to\R[\Lp]_m$. 

Note that $m+r=\dim\Lp$; in the sequel, it will be convenient to identify 
the representations $\Lambda^m\Lp\cong\Lambda^r\Lp$ and view 
$\Phi:\Lambda^r\Lp\to\R[\Lp]_m$ 
(this can be done because $\Lp$ is a real orthogonal 
representation of $K$). Fix a basis $\{y_1,\ldots,y_d\}$ of $\Lp$. 
Inspired by~\cite{gru}, we
define a linear map $A:\Lambda^r\Lp\to \Lk\otimes\Lambda^{r-2}\Lp$ by setting
\[ A(y_{i_1}\wedge\cdots\wedge y_{i_r}) = \sum_{a<b}[y_{i_a},y_{i_b}]
\otimes y_{i_1}\wedge\cdots\hat{y}_{i_a}\wedge\cdots\wedge\hat{y}_{i_b}
\wedge\cdots\wedge y_{i_r}. \]
Then $A$ is equivariant and one easily computes
its adjoint $A^*:\Lk\otimes\Lambda^{r-2}\Lp\to\Lambda^r\Lp$ to be 
\[ A^*(x\otimes y_{i_1}\wedge\cdots\wedge y_{i_{r-2}})
=\sum_{\ontop{a<b}{a,b\neq i_1,\ldots,i_{r-2}}}
\langle \mathrm{ad}_xy_a,y_b\rangle
y_{i_1}\wedge\cdots\wedge y_a\wedge\cdots\wedge y_b
\wedge\cdots\wedge y_{i_{r-2}}. \]

\begin{prop}\label{Phi-A*}
The composite map $\Phi\circ A^* =0$.
\end{prop}

\Pf Fix a Cartan subspace $\La\subset\Lp$ and an orthonormal basis
$\{y_i\}_{i=1}^d$ of $\Lp$ such that $y_1,\ldots,y_r\in\La$.
Since the $K$-orbits meet $\La$ orthogonally, for $a\in\La$ we have
\[ \Phi(y_{i_1}\wedge\cdots y_a\wedge\cdots\wedge y_b
\wedge\cdots\wedge y_{i_r})(a)\neq 0 \]
if and only if $y_{i_1}\wedge\cdots y_a\wedge\cdots\wedge y_b
\wedge\cdots\wedge y_{i_r}=y_1\wedge\cdots\wedge y_r$. 
However, in this case, $y_a$, $y_b\in\mathcal\La$ implies
that $\langle \mathrm{ad}_{x}y_a,y_b\rangle=0$ for every $x\in\Lk$. 
This shows that $\Phi\circ A^*|_{\mathfrak a}=0$. 
Let $\{\alpha_i\}$ be an orthonormal basis of 
$\Lk\otimes\Lambda^{r-2}\Lp$. Since $\Phi\circ A^*$ is 
$K$-equivariant, Lemma~\ref{inv} says that 
$f=\sum_i[(\Phi\circ A^*)(\alpha_i)]^2$ is $K$-invariant, and we already know
that $f|_{\mathfrak a}=0$. It follows that $f=0$ on $\Lp$. Hence 
$\Phi\circ A^*(\alpha_i)=0$ for all $i$ and thus $\Phi\circ A^*=0$. \EPf

\medskip

Since $(\mathrm{im}A^*)^\perp=\ker A$, in view of Proposition~\ref{Phi-A*}
we need only to consider the restriction of $\Phi$ to the kernel of $A$.
It is often easier to deal with complex representations,
so we complexify everything. Now our problem is equivalent 
to identifying $\Lk^c$-irreducible components of 
$\ker A^c\subset\Lambda^r(\Lp^c)$ which are not mapped to zero under $\Phi^c$. 
Let $\Lt$ be the Lie algebra of a maximal torus of $K$.
Then we have the (complex) root space decomposition
\[ \Lk^c = \Lt^c + \sum_{\alpha\in\Delta_{\mathfrak {k^c}}}(\Lk^c)_\alpha \]
and the weight space decomposition
\[ \Lp^c = (\Lp^c)_0 + \sum_{\alpha\in\Delta_{\mathfrak p^c}}(\Lp^c)_\alpha, \]
where $(\Lp^c)_0$ is the centralizer of $\Lt^c$ in $\Lp^ c$,
$\Lt^c+(\Lp^c)_0$ is a Cartan subalgebra of $\Lg^c$ and 
$\dim(\Lk^c)_\alpha=\dim(\Lp^c)_\alpha=1$ (see e.g.~\cite[p.15]{pan}).
Choose a system of positive roots $\Delta_{\mathfrak {k^c}}^+
\subset\Delta_{\mathfrak {k^c}}$ , and let 
$\Lb=\Lt^c+\sum_{\alpha\in\Delta_{\mathfrak {k^c}}^+}(\Lk^c)_\alpha $ 
denote the corresponding Borel subalgebra of $\Lk^c$. 
It is clear that the map $\Lv\mapsto[\Lambda^r\Lv]$ sets
up a bijective correspondence between 
the $r$-dimensional $\Lb$-stable subspaces of $\Lp^c$ and 
the highest weight vectors of $\Lambda^r(\Lp^c)$. 

One often combines the aforementioned Chevalley restriction theorem 
with the Chevalley theorem for finite reflection groups, which says that  
the algebra of invariants $\R[\La]^{W(\mathfrak a)}$ is a free polynomial
algebra, namely it has $r$ algebraicaly independent homogeneous generators 
$f_1,\ldots,f_r$. Now $f_1,\ldots,f_r$ are also
algebraicaly independent homogeneous generators for
$\C[\Lp^c]^{K^c}$.
By a result of Panyushev~\cite{pan2}, 
$y\in\Lp^c$ is regular if and only if the set of linear forms on $\Lp^c$
$\{(df_1)_y,\ldots,(df_r)_y\}$ is linearly independent. In this case,
one can find a $r$-dimensional complex subspace $\Lv\subset\Lp^c$ 
(necessarily transversal to the orbit $K^c\cdot y$) such that 
the set of restricted linear forms $\{(df_1)_y|_{\mathfrak v},\ldots,
(df_r)_y|_{\mathfrak v}\}$ is linearly independent. 
Since $\Lv$ is transversal to $K^c\cdot y$, $\Phi^c([\Lambda^r\Lv])\neq0$. 
This remark is effective in the case $\Lv$ can be taken to be 
a $\Lb$-stable subspace of $\Lp^c$, for in that case 
$[\Lambda^r\Lv]$ is a highest weight vector of $\Lambda^r(\Lp^c)$,
whence determines an irreducible component not mapped to 
zero under $\Phi^c$. We summarize this discussion in the 
the following theorem.  

\begin{thm}\label{alg}
Let $\Lv$ be a $\Lb$-stable subspace of $\Lp^c$ 
such that for some $y\in\Lp^c$ 
the set of linear forms on $\Lp^c$ $\{(df_1)_y|_{\mathfrak v},\ldots,
(df_r)_y|_{\mathfrak v}\}$ is linearly independent. Then 
$[\Lambda^r\Lv]$ is the highest weight vector of an irreducible
component of $\Lambda^r(\Lp^c)$ which is not mapped to zero under 
$\Phi^c$.
\end{thm}

\begin{rmk}\label{rmk}
\em
In this remark, assume for simplicity 
that $(K,\Lp)$ is irreducible. Then the Casimir 
element $\omega$ of $\Lk^c$ with respect to the Cartan-Killing form of 
$\Lg^c$ can be normalized so as to act as the identity on $\Lp^c$ 
(even if $(K,\Lp)$ is not absolutely irreducible). 
We quote a result of Kostant~\cite{kos} and Panyushev~\cite{pan} stating that:
\emph{The maximal eigenvalue of $\omega$ on $\Lambda^r(\Lp^c)$ is $r$. 
The corresponding
eigenspace $\mathcal M$ is spanned by decomposable elements. Indeed,
$\mathcal M$ is spanned by $[\Lambda^r\mathfrak d]$, where $\mathfrak d$ 
runs over
all $r$-dimensional Abelian subalgebras of $\Lp^c$.}
Hence $\mathcal M\subset\ker A$, but in general the inclusion is strict,
as there are simple examples with $r=3$ in which $\ker A$ is not spanned 
by decomposable elements. In any case, $\mathcal M$ is a first approximation
to $\ker A$. 
\end{rmk}

\section{Applications}

In this section we apply Theorem~\ref{alg} to effectively
compute upper bounds for the minimum number of squares 
required to represent $\delta$ in some concrete examples.
In particular, we recover Domokos estimate in the case of the 
discriminant of $n\times n$ real symmetric matrices~\cite[Thm.~6.2]{dom}. 

We start with the symmetric space $G/K=Sp(n,\R)/U(n)$. Its isotropy
representation is the realification of the representation of $U(n)$ on 
$V=\mathrm{Sym}(n,\C)$ given by $\rho(g)X=gXg^t$, 
so the associated discriminant
can also be interpreted as the discriminant of 
complex $n$-ary quadratic forms. As Cartan subspace, one can take 
the set of real diagonal matrices. 

The complexified symmetric space
$G^c/K^c$ has~\cite[p.594]{gw} 
\[ G^c=Sp(2n,\C)=\{g\in M(n,\C):g^t\,Jg=J\}, \]
the complex symplectic group with respect to
$J=\left(\begin{array}{cc}0&I_n\\-I_n&0\end{array}\right)$
where $I_n$ denotes an $n\times n$ identity block. As involution of
$G^c$, we have conjugation by
$\left(\begin{array}{cc}I_n&0\\0&-I_n\end{array}\right)$,
which yields that $K^c\cong GL(n,\C)$ consisting of the
matrices 
\[ \left(\begin{array}{cc}g&0\\0&(g^t)^{-1}\end{array}
\right)\qquad\mbox{with}\qquad g\in GL(n,\C), \]
whereas $\Lp^c$ consists of the matrices ($n\times n$ blocks)
\begin{equation}\label{y}
 y=\left(\begin{array}{cc}0&X\\Y&0\end{array}
\right)\qquad\mbox{with}\qquad X^t=X,\ Y^t=Y. 
\end{equation}
Further, $(K^c,\Lp^c)$ is equivalent to $\rho\oplus\rho^*$;
here we can identify $V^*$ with $V$
as a vector space and then~$\rho^*(g)X=(g^t)^{-1}Xg^{-1}$.
Note that $\rho$ is irreducible with highest weight $2\theta_1$.
  
The polynomials $f_j(y)=\mathrm{tr}((XY)^j)$ for $j=1,\ldots,n$ and 
$\nu$ as in~(\ref{y}) form a complete
set of invariants. One easily computes
\[ (df_j)_y(\tilde X,\tilde Y)=j\left(\langle
(XY)^{j-1}X,\tilde Y\rangle+\langle Y(XY)^{j-1},\tilde X\rangle\right), \]
where $(\tilde X,\tilde Y)\in V\oplus V^*$ and $\langle\cdot,\cdot\rangle$
denotes a constant multiple of 
the Cartan-Killing form of $\Lg^c$. In particular,
by taking $Y$ to be the identity and $\tilde Y=0$, we get 
$(df_j)_y(\tilde X,0)=j\,\langle X^{j-1},\tilde X\rangle$, whence
$y$ is regular if and only if $\{I,X,\ldots,X^{n-1}\}$ is a linearly
independent set of matrices. In the following, we construct $X$ 
such that already the first columns of the matrices $I,X,\ldots,X^{n-1}$
form a linearly independent set. It follows that the corresponding
$y$ is regular and moreover we can take the subspace 
$\Lv$ of $\Lp^c$ as in Theorem~\ref{alg}
as being the subspace of $V$ formed by symmetric matrices $\tilde X$ whose
only nonzero entries lie in the first column or first line. 

Let $Z$ be a real symmetric matrix with distinct eigenvalues,
say diagonal. It is known that such a matrix 
admits a cyclic vector, 
that is a vector $w_1\in\R^n$ such that 
$\{w_1,Zw_1,\ldots,Z^{n-1}w_1\}$ is a basis of $\R^n$.
We can assume $w_1$ is a unit vector and then complete it  
to form an orthonormal basis $\{w_1,\ldots,w_n\}$ of $\R^n$.  
Set $X=M^{-1}ZM$, where $M$ is the orthogonal matrix whose column
vectors are the $w_j$. Then $X^{j-1}$ is a symmetric matrix
whose first column is
\[ \left(\begin{array}{c}w_1^t\,Z^{j-1}w_1\\ \vdots \\
w_n^t\,Z^{j-1}w_1 \end{array}\right), \]
namely, contains the coordinates of the vector $Z^{j-1}w_1$
with respect to the basis  $\{w_1,\ldots,w_n\}$; hence the set 
of such columns for $j=1,\ldots, n$ form a linearly independent
set, as we wished.

We have that $E_{ij}+E_{ji}\in V$ is a weight vector of 
$\rho$ of weight $\theta_i+\theta_j$.
It follows that
$[\Lambda^n\Lv]$ is a highest weight vector of $\mathfrak{gl}(n,\C)$
with highest weight (compare~\cite[p.274]{gw})
\[ 2\theta_1+(\theta_1+\theta_2)+\cdots+(\theta_1+\theta_n)=n\theta_1
+\underbrace{(\theta_1+\cdots+\theta_n)}_{\mbox{$=0$ on 
$\mathfrak{sl}(n,\C)$}}, \]
so the corresponding representation of $U(n)$ is the $n$-th
symmetric power $S^n(\C^n)$ of the standard representation on $\C^n$. 
Since the real dimension of this representation is $2\,\binom{2n-1}{n}$, 
this proves Theorem~\ref{spn-un} stated in the introdution. 

The lowest nontrivial case $n=2$ is special, for 
the theorem then says $\delta$ is the sum of 
six squares; however, it is easy to see the 
decomposition into $U(2)$-irreducible representations
$\Lambda^2(S^2(\C^2)^r)\cong(\C^3)^r\oplus\R^5\oplus\R^3\oplus\R$,
where $\Phi$ kills exactly the last two components (corresponding
to the adjoint action): we have that $(\C^3)^r=(S^2(\C^2))^r$ yields
$\delta$ as a sum of six squares as in 
Theorem~\ref{spn-un}, but $\R^5$, a real form of $S^4(\C^2)$ 
on which the center of $U(2)$ does not act, yields the better 
result that $\delta$ is a sum of five squares. 
Moreover,
one computes directly for a 
complex symmetric matrix $\left(\begin{array}{cc}z_1&z_3\\z_3&z_2
\end{array}\right)$ that  
\[ \delta = |z_1z_2-z_3^2|^2\left[\,\left(|z_1|^2-|z_2|^2\right)^2
+4\left|z_1\bar z_3
+\bar z_2z_3\right|^2\,\right], \]
which shows in fact the best result that 
$\delta$ is a sum of two squares. 
For $n\geq3$, Theorem~\ref{spn-un} is probably neither optimal,
though it should get closer to that. 

Finally, we quickly revisit the case of $SO(n)$-conjugation of
traceless real symmetric matrices,
or, the isotropy representation of $SL(n,\R)/SO(n)$.
It is usual and convenient to model the 
complexified representation $(SO(n,\C),\mathrm{Sym}_0(n,\C))$ using the
nondegenerate symmetric bilinear form given by 
\[ Q=\left(\begin{array}{cc}0&I_\ell\\I_\ell&0\end{array}\right)
\quad\mbox{if $n=2\ell$}\qquad\mbox{and}\qquad 
Q=\left(\begin{array}{ccc}0&I_\ell&0\\I_\ell&0&0\\0&0&1\end{array}\right)
\quad\mbox{if $n=2\ell+1$,} \]
so that $SO(n,\C)\cong SO(Q):=\{g|g^t\,Qg=Q\}$ 
with Lie algebra $\mathfrak{so}(n,\C)\cong
\mathfrak{so}(Q):=\{A|A^t\,Q+QA=0\}$ and 
$\mathrm{Sym}_0(n,\C)\cong V:=\{X|X^t=QXQ^{-1},\ \mathrm{tr} X=0\}$. 
Now $(K^c,\Lp^c)$ is the 
representation $\rho$ of $SO(Q)$ on $V$
given by $\rho(g)X=gXg^{-1}$. We have that 
$F:S^2_0(\C^n)\to V$ given by 
$F(v\cdot w)=\frac12(vw^t+wv^t)Q$ is an equivariant isomorphism;
here $S^2_0(\C^n)$ denotes the nontrivial component of the 
symmetric square $S^2(\C^n)$. 
It follows that $\rho$ is irreducible
with highest weight $2\theta_1$. 

The polynomials $f_j(X)=\mathrm{tr}(X^j)$ for $j=2,\ldots,n$ form
a complete set of invariants. One easily computes that 
\[ (df_j)_X(\tilde X)=j\langle X^{j-1},\tilde X\rangle, \]
where $\tilde X\in V$ and $\langle\cdot,\cdot\rangle$
denotes a constant multiple of the Cartan-Killing form of $\Lg^c$. Therefore 
$X$ is regular if and only if $\{X,X^2,\ldots,X^{n-1}\}$ is a 
linearly independent set of matrices. 

The matrix $X$ chosen in~\cite[p.13]{dom} (denoted $A$ there) 
acting on 
the canonical basis $\{e_1,\ldots,e_n\}$ of $\C^n$ as
\[ e_1\mapsto e_{\ell+1}\mapsto e_{\ell+2}\mapsto\cdots\mapsto
e_n\mapsto e_\ell\mapsto e_{\ell-1}\mapsto\cdots\mapsto e_2\mapsto e_1 \]
lies in $V$ and has $e_1$ as cyclic vector.
Thus the first columns of the matrices 
$I,X,\ldots,X^{n-1}$ form a linearly independent set.
Since $X^{n-1}=X^{-1}=X^t$, the set 
$\{I,X,\ldots,X^{n-1}\}$ is invariant under transposition
and hence
also the first lines of 
$I,X,\ldots,X^{n-1}$ form a linearly independent set.
In particular, taking out the identity $I$ from this set,
it suffices to consider first lines minus the $(1,1)$-entry of the
matrices $X,\ldots,X^{n-1}$ to see that they
form a linearly independent set.
Thus $X$ is regular and we can take the subspace 
$\Lv$ of $\Lp^c$ as in Theorem~\ref{alg}
as being the subspace of $V$ formed by matrices $\tilde X$ whose
only nonzero entries are in positions $(1,2)$ through $(1,n)$. 
Note that  $\Lv=F(\langle e_1^2,e_1e_2,\ldots,\widehat{e_1e_{\ell+1}},
\ldots,e_1e_n\rangle)$, and the weight of 
$e_i$ in $\C^n$ is 
\[ \begin{array}{ll}
\theta_i&\mbox{if $i=1,\ldots,\ell$,}\\
-\theta_i&\mbox{if $i=\ell+1,\ldots,2\ell$,}\\
0&\mbox{if $i=2\ell+1$.}\\
\end{array} \]
It follows that 
$[\Lambda^{n-1}\Lv]$ is a highest weight vector of weight 
\[ \begin{array}{ll}
 2\theta_1+(\theta_1+\theta_2)+\ldots+(\theta_1+\theta_\ell)+ 
(\theta_1-\theta_2)+\ldots+(\theta_1-\theta_\ell)&\mbox{if $n=2\ell$,}\\
2\theta_1+(\theta_1+\theta_2)+\ldots+(\theta_1+\theta_\ell)+ 
(\theta_1-\theta_2)+\ldots+(\theta_1-\theta_\ell)+\theta_1
&\mbox{if $n=2\ell+1$,}
\end{array} \]
which in both cases equals $n\theta_1$, which corresponds to the 
representation on the 
space of $n$-variable spherical harmonics of 
degree~$n$, of dimension
$\dim S^n(\R^n)-\dim S^{n-2}(\R^n)$.
This proves Theorem~\ref{dom} stated in the introduction.

\providecommand{\bysame}{\leavevmode\hbox to3em{\hrulefill}\thinspace}
\providecommand{\MR}{\relax\ifhmode\unskip\space\fi MR }
\providecommand{\MRhref}[2]{%
  \href{http://www.ams.org/mathscinet-getitem?mr=#1}{#2}
}
\providecommand{\href}[2]{#2}

\end{document}